\UseRawInputEncoding
\documentclass[12pt,reqno]{amsart}
\usepackage{amssymb,amsmath,graphicx,
	amsfonts,euscript}
\usepackage{color}
\usepackage{mathrsfs}
\theoremstyle{plain}

\newtheorem{theorem}{Theorem}[section]
\newtheorem{lemma}[theorem]{Lemma}

\newtheorem{remark}[theorem]{Remark}
\numberwithin{equation}{section}

\begin{document}

\title[The symmetry of steady stratified periodic water waves]{ The symmetry for two class of steady stratified periodic water waves}

\author[Fengquan Li, Fei Xu, Yong Zhang]{Fengquan Li, Xu Fei, Zhang Yong}

\address[Fengquan Li]{School of Mathematical Sciences, Dalian University
of Technology, Dalian, 116024, China}

\address[Fei Xu]{School of Mathematical Sciences, Dalian University
of Technology, Dalian, 116024, China}

\address[Yong Zhang]{ School of Mathematical Sciences, Dalian University
of Technology, Dalian, 116024, China} \email{18842629891@163.com}

\begin{abstract}
In this paper, we mainly consider two class of travelling stratified periodic water waves, one with negative (or without) surface tension and the other with constant Bernoulli's function and stagnation points. We first establish the symmetry result for stratified water waves with negative (or without) surface tension, but without stagnation by using the modified maximum principle. Furthermore, the symmetry property of stratified water waves with constant Bernoulli's function and stagnation points is also obtained provided the monotonic property is known.
	\end{abstract}
\maketitle
\noindent {\bf Key Words:}
{Symmetry, Stratified water wave, Negative surface tension, Stagnation points}

\section{\bf Introduction}
The density stratification in the oceanography and geophysical fluid is ubiquitous due to the change of water temperature and salinity at different depth. The phenomenon of stratification has inspired many scholars to investigate its mechanism. The paper mainly considers the qualitative properties of two class of steady stratified periodic water waves.

In fact, the symmetry property of water waves has been a complex topic, especially for stratified waves. It was shown in \cite{17} that irrotational gravity water waves with profile having a unique minimum and maximum within a period were symmetric, where the maximum principles for elliptic problems and the moving plane method (see \cite{21}) played a key role. Moreover, the symmetry of irrotational water waves with negative surface tension was also discussed by Okamoto and Sh\={o}ji in \cite{17}. When considering the periodic gravity water waves with vorticity, the symmetry property was extended by Constantin and Escher to finite depth \cite{15} and deep water \cite{18} flows without stagnation points under some restriction on the derivative of the vorticity function. Later, Hur \cite{24} established the symmetry for an arbitrary vorticity but requires a quite precise knowledge of all the streamlines in the fluid. By introducing hodograph transformation, the result was further improved in \cite{5} for any arbitrary vorticity functions. Following the breakthrough work in \cite{5}, the symmetry properties of periodic water waves have been extended to more general rotational flows, such as flows are stratified \cite{3,6}, or with negative surface tension \cite{19} and flows are near equatorial region \cite{20,29}, or with arbitrarily many crests and troughs \cite{28}. Besides, the moving plane method was also applied to investigate the symmetry properties of gravity solitary waves in \cite{22,23}.

In contrast to most previous studies on symmetry, in this paper we admit that there exist negative surface tension or internal stagnation points but not overhanging profiles in steady stratified periodic water waves. As for the existence of these waves, we recommend the work \cite{11,12,25,26,27,14}. These characteristics would bring some essential difficulties to study the symmetry property. It is not straightforward to extend the symmetry approach presented in \cite{5} to stratified waves since one obtains elliptic operators with a zero-order term of the bad sign. The author \cite{3} discussed the symmetry of stratified periodic by restricting the amplitude of water waves to use the classical maximum principle. As Constantin mentioned in his book \cite{7} "it seems that some restriction on the size of the free term cannot be avoided to apply maximum principles". However, here we establish the symmetry of stratified water waves by introducing some modified maximum principles without any restrictions on the elevation of stratified water waves, which can be extended to stratified waves with negative surface tension. In addition, when there exist internal stagnation points in stratified periodic water waves, the symmetry result is also obtained by developing the idea in \cite{15}.

The rest of this paper is arranged as follows. In Sect.2, we recall the setup of governing equations for two-dimensional steady stratified periodic water waves. In Sect.3, we use the modified maximum principle, which is suitable for stratified water waves of height function formulation without stagnation points, to establish the symmetry of water waves. In Sect.4, we first investigate the existence and location of stagnation points of the steady stratified periodic water waves with constant Benoulli's function and then obtain the symmetry property based on the pseudo stream function formulation.

\section{\bf Formulation of the problem}

Fix a Cartesian coordinate system so that the $X$-axis points to be horizontal, and the $Y$-axis to be vertical. Assume that the floor of the sea bed is flat and occurs at $Y=0$, and $Y=\eta(t,X)>0$ be the free surface. We shall normalize $\eta$ by choosing the axis so that the free surface is oscillating around the line $Y=d>0$. As usual we let $u=u(t,X,Y)$ and  $v=v(t,X,Y)$ denote the horizontal and vertical velocities respectively, and let $\rho=\rho(t,X,Y)> 0$ be the density and $P=P(t,X,Y)$ be the pressure, all of which have the form $(X-ct)$ in steady travelling waves, where $c$ represents the speed of wave. For convenience, we denote $x=X-ct, y=Y$ and consider the problem in $\Omega=\{(x,y)|-\pi<x<\pi, 0<y<\eta(x)\}$ due to the periodicity.

For water waves, the incompressibility is represented mathematically by the requirement that the vector field $(u,v)$ is divergence free for all time
\begin{eqnarray}
u_{X}+v_{Y}=0. \label{eq2.1}
\end{eqnarray}
The conservation of mass implies
\begin{eqnarray}
\rho_{t}+(\rho u)_{X}+(\rho v)_{Y}=0. \label{eq2.2}
\end{eqnarray}
Then the relation of time-space $(X-ct)$ and (\ref{eq2.1})-(\ref{eq2.2}) demonstrate
\begin{eqnarray}
(u-c)\rho_{x}+v\rho_{y}=0. \label{eq2.3}
\end{eqnarray}
Therefore, the governing equations in velocity field formulation are expressed by the nonlinear free-boundary problem (see\cite{1})
\begin{eqnarray}
\left\{\begin{array}{llll}{u_{x}+v_{y}=0} & {\text { for }0 \leq y \leq \eta_(x)}, \\
{(u-c)\rho_{x}+v\rho_{y}=0} & {\text { for } 0 \leq y \leq \eta_(x)}, \\
{\rho(u-c) u_{x}+ \rho v u_{y}=-P_{x}} & {\text { for }0 \leq y\leq \eta(x)}, \\
{\rho(u-c) v_{x}+\rho v v_{y}=-P_{y}-g\rho} & { \text { for }0 \leq y \leq \eta(x)}, \\
{v=0} & {\text { for } y=0}, \\
{v=(u-c) \eta_{x}} & {\text { on } y=\eta(x)}, \\
{P=P_{atm}-\sigma\frac{\eta_{xx}}{(1+\eta_{x}^{2})^{\frac{3}{2}}}} & {\text { on } y=\eta(x)},\end{array}\right. \label{eq2.4}
\end{eqnarray}
where $P_{atm}$ is the constant atmosphere pressure, $g=9.8 m/s^{2}$ is the (constant) gravitational acceleration at the Earth's surface and the parameter $\sigma$ is the coefficient of surface tension (here we choose it's not positive).

Observe that, by conservation of mass and incompressibility, $\rho$ is transported and the vector field is divengence free. Therefore we may introduce a (relative) pseudo-stream function $\psi=\psi(x,y)$ satisfying
\begin{eqnarray}
\psi_{x}=-\sqrt {\rho}v, ~~~~~~~\psi_{y}=\sqrt {\rho}(u-c). \label{eq2.5}
\end{eqnarray}

It is a straightforward calculation to check that $\psi$ is indeed a (relative) stream function in the usual sense, i.e. its gradient is orthogonal to the vector field in the moving frame at each point in the fluid domain. As usual, we shall refer to the level sets of $\psi$ as the streamlines of the flow.
Define the relative pseudo mass flux (relative to the uniform flow at speed $c$) as
\begin{eqnarray}
p_{0}=\psi(x,\eta(x))-\psi(x,0)= \int_{0}^{\eta(x)}\sqrt {\rho(x,y)}[u(x,y)-c]dy,  \label{eq2.6}
\end{eqnarray}
Then boundary conditions in (\ref{eq2.4}) imply that
\begin{eqnarray}
\begin{aligned}
\frac{dp_{0}}{dx}
&=\eta_{x}(\sqrt{\rho(x,\eta(x))}[u(x,\eta(x))-c])+\int_{0}^{\eta(x)}\partial_{x}(\sqrt {\rho(x,y)}[u(x,y)-c])dy \nonumber\\
&=\sqrt{\rho(x,\eta(x))}v(x,\eta(x))-\int_{0}^{\eta(x)}\partial_{y}(\sqrt {\rho(x,y)}v)dy=0,
\end{aligned}
\end{eqnarray}
which indicates that $p_{0}, \psi(x,\eta(x))$ and $\psi(x,0)$ are constants. Without loss of generality, we choose $\psi\equiv 0$ on the free boundary $y=\eta(x)$, which forces $\psi= -p_{0}$ on $y=0$. Since $\rho$ is transported, it must be constant on the streamlines. We may therefore let streamline density function $\rho\in C^{1,\alpha}([p_{0}, 0]; R^{+})$ be given by
\begin{eqnarray}
\rho(x,y)=\rho(-\psi(x,y)). \label{eq2.7}
\end{eqnarray}
throughout the fluid.
From Bernoulli's law, we know that
\begin{eqnarray}
E=\frac{\rho}{2}((u-c)^{2}+v^{2})+gy\rho+P \label{eq2.8}
\end{eqnarray}
is a constant along each streamline. Then there exists a function $\beta \in C^{1,\alpha}([0, |p_{0}|]; R)$ such that
\begin{eqnarray}
\frac{dE}{d\psi}=-\beta(\psi),  \label{eq2.9}
\end{eqnarray}
where $\beta$ is called the Bernoulli function corresponding to the flow (see \cite{1}). Physically it describes the variation of specific energy as a function of the streamlines. It is worth noting that when $\rho$ is a constant, $\beta$ reduces to the vorticity function. It's obvious that (\ref{eq2.8}) and (\ref{eq2.9}) show
\begin{eqnarray}
\frac{dE}{d\psi}=\Delta\psi-gy\rho'(-\psi)=-\beta(\psi(x,y)). \label{eq2.10}
\end{eqnarray}
Moreover, evaluating Bernoulli's theorem on the free surface, we find
\begin{eqnarray}
2E|_{\eta}=2P_{atm}+|\nabla\psi|^{2}+2g\eta\rho(-\psi)-2\sigma\frac{\eta_{xx}}{(1+\eta_{x}^{2})^{\frac{3}{2}}}~~~~~~on ~~y=\eta(x). \label{eq2.11}
\end{eqnarray}

Summarizing the above considerations, we can reformulate the governing equations (\ref{eq2.4}) as follows
\begin{eqnarray}
\left\{\begin{array}{ll}{ \Delta\psi-gy\rho'(-\psi)=-\beta(\psi)} & {\text { in }0 < y<\eta(x)}, \\
{|\nabla\psi|^{2}+2g\rho(-\psi)y-2\sigma\frac{\eta_{xx}}{(1+\eta_{x}^{2})^{\frac{3}{2}}}=Q} & { \text { on } y=\eta(x)}, \\
{\psi=0} & {\text { on } y=\eta(x)}, \\
{\psi=-p_{0}} & {\text { on } y=0},\end{array}\right. \label{eq2.12}
\end{eqnarray}
where $Q=2(E|_{\eta}-P_{atm})$ is a constant.

\section{\bf The stratified water waves with $\psi_{y}<0$ and $\sigma\leq0$}
In this section, we are concerned with the steady stratified periodic water waves with negative (or without) surface tension. From the physical view, maybe considering the positive surface tension is more meaningful. It is a pity that establishing the symmetry in this case is still an open problem even for steady periodic water waves without stratification (see \cite{7}). While the following results are more mathematically interesting. Considering the stability of flows, we assume that the streamline density function $\rho$ is nonincreasing (i.e $\rho'(-\psi)\leq 0$). Besides, we also require that the horizontal current velocity is smaller than wave speed (i.e $\psi_{y}<0$), which admits to introduce the Dubreil-Jacotin's transformation \cite{13} by
$$
q=x,~~~~~p=-\psi(x,y).
$$
This transforms the fluid domain
$$\Omega=\{(x,y):~x\in(-\pi,\pi),~0<y<\eta(x)\}$$
into rectangular domain
$$D=\{(q,p):~-\pi< q < \pi,~ p_{0}< p< 0\},$$
where $p_{0}<0$ and let $\overline{D}$ denote its closure.
Moreover, $\beta$ and $\rho$ in (\ref{eq2.7}) and (\ref{eq2.9}) can be written as
$$
\beta=\beta(-p), ~~~~~~\rho=\rho(p).
$$
Defining the height above the flat bed by
$$
h(q,p):=y
$$
and taking the mean of $h(q,p)$ along the free surface, we get
$$
d=\int_{-\pi}^{\pi}h(q,0)dq.
$$
In addition, we can deduce
$$
\psi_{y}=-\frac{1}{h_{p}},~~~~ \psi_{x}=\frac{h_{q}}{h_{p}},
$$
$$
\partial_{q}=\partial_{x}+h_{q}\partial_{y},~~~~ \partial_{p}=h_{p}\partial_{y}.
$$
Consequently, we can rewrite the governing equations (\ref{eq2.12}) as the height function formulation
\begin{eqnarray}
\left\{\begin{array}{lll}{\left(1+h_{q}^{2}\right) h_{p p}-2 h_{q} h_{p} h_{q p}+h_{p}^{2} h_{q q}+[\beta-g(h-d)\rho']h_{p}^{3}=0}  & {\text { in } p_{0}<p<0}, \\
{1+h_{q}^{2}+(2g\rho h-Q-2\sigma\frac{h_{qq}}{(1+h_{q}^{2})^{\frac{3}{2}}})h_{p}^{2}=0} & {\text { on } p=0}, \\
{h=0} & {\text { on } p=p_{0}}.\end{array}\right. \label{eq3.1}
\end{eqnarray}
Based on the formulation (\ref{eq3.1}), we can obtain the following symmetry result.
\begin{theorem}\label{thm3.1}
The steady stratified periodic water waves with negative (or without) surface tension are symmetric if streamlines are monotonic between the trough line and the crest line and they are strictly monotonic in a neighborhood of the trough line $x=\pm\pi$.
\end{theorem}

\begin{remark}:
The existence of stratified water waves without surface tension has been rigorously proved in \cite{1}. We can refer to \cite{11,12} (up to a slight modification) for the existence of steady stratified periodic water waves with negative surface tension (see \cite{19} for the gravity water waves with negative surface tension by modifying the result in \cite{30,31}).
\end{remark}
Before starting the proof, let's introduce a modified maximum principle.
\begin{lemma} \label{lem3.1}(cf. \cite{2})
 Suppose $D\subset R^{2}$ be an open rectangle and $\omega\in C^{2}(D)\cap C(\overline{D})$ satisfy $L\omega\leq0$ for the uniformly elliptic operator $L=a_{ij}\partial_{ij}+b_{i}\partial_{i}+c$ with $a_{ij},b_{i},c\in C(\overline{D})$. If $\omega\geq 0$ in $D$, then the followings hold \\
(1) The weak maximum principle: $\omega$ attains its minimum on $\partial D$. \\
(2) The strong maximum principle: If $\omega$ attains its minimum in $D$, then $\omega$ is constant in $\overline{D}$.\\
(3) Hopf's principle: Let Q be a point on $\partial D$, different from the corners of the rectangle $\overline{D}$. If $\omega(Q)<\omega(X)$ for all $X$ in $D$, then $\partial_{\vec{\nu}}\omega(Q)< 0$ where $\nu$ is the outward normal direction at $Q$.
\end{lemma}
\begin{remark}:
In general, Lemma \ref{lem3.1} holds for the uniformly elliptic operator $L$ with a zeroth order term $c(x)\leq 0$, however, the additional condition $\omega\geq0$ in $D$ would remove the restriction for the sign of $c(x)$. Indeed, let $c(x)=c^{+}(x)-c^{-}(x)$, where $c^{+}:=\max\{0,c\}$ and $c^{-}:=\max\{0,-c\}$, then $\tilde{L}\omega:=a_{ij}\partial_{ij}\omega+b_{i}\partial_{i}\omega-c^{-}\omega=L\omega-c^{+}w\leq0$ in $D$ provided $L\omega\leq0, \omega\geq0$ in $D$.
\end{remark}
~\\
{\bf Proof of the Theorem \ref{thm3.1}:}
Let $h,\tilde{h}\in C^{2,\alpha}(\bar{D})$ be two solutions to the water wave problem (\ref{eq2.23}) without the stagnation points (i.e. $h_{p}=\frac{1}{\sqrt{\rho}(c-u)}>0$), then the difference is a solution of
$$
L(h-\tilde{h})=0,~~~in ~~~D,
$$
where
\begin{align}\label{eq3.2}
 L
 &=(1+h_{q}^{2})\partial_{pp}-2h_{p}h_{q}\partial_{qp}+h_{p}^{2}\partial_{qq} \nonumber\\
 &+[\tilde{h}_{qq}(h_{p}+\tilde{h}_{p})-2\tilde{h}_{q}\tilde{h}_{qp}+(\beta-g\rho_{p}(\tilde{h}-d))(h_{p}^{2}+h_{p}\tilde{h}_{p}+\tilde{h}_{p}^{2})]\partial_{p} \nonumber\\
 &+[\tilde{h}_{pp}(h_{q}+\tilde{h}_{q})-2h_{p}\tilde{h}_{qp}]\partial_{q}-g\rho_{p}h_{p}^{3}
 \end{align}
is a uniformly elliptic operator with continuous coefficients. Without loss of generality, let's assume that the wave trough is at $q=\pm \pi$. From the construction of solutions, it's clear that, for any $\lambda$, the function $\tilde{h}(q,p)=h(2\lambda-q,p)$ also satisfies the system (\ref{eq3.1}) if $h(q,p)\in C^{2,\alpha}_{per}(\bar{D})$ with $h_{p}>0$ is a solution to (\ref{eq3.1}), where $2\lambda-q$ is called the reflection of $q$ about $\lambda$.
 Define the associated reflection function by
$$
w(q,p;\lambda)=h(q,p)-\tilde{h}(q,p)=h(q,p)-h(2\lambda-q,p)
$$
and it's easy to check that $w$ satisfies
\begin{eqnarray}\label{eq3.3}
\left\{\begin{array}{lll}{Lw=0}  & {\text { in } p_{0}<p<0}, \\
{h_{q}^{2}-\tilde{h}_{q}^{2}+(2g\rho h-Q-2\sigma\frac{h_{qq}}{(1+h_{q}^{2})^{\frac{3}{2}}})h_{p}^{2}-(2g\rho \tilde{h}-Q-2\sigma\frac{\tilde{h}_{qq}}{(1+\tilde{h}_{q}^{2})^{\frac{3}{2}}})\tilde{h}_{p}^{2}=0}& {\text { on } p=0}, \\
{w=0} & {\text { on } p=p_{0}}.\end{array}\right. \label{eq2.23}
\end{eqnarray}

To finish the proof of the Theorem \ref{thm3.1}, it suffices for us to show that
$$
w(q,p;0)\equiv 0~~~for~~~(q,p)\in[0,\pi]\times[p_{0},0].
$$
For a clear presentation, let's divide the process into following two steps.\\
~\\
{\bf Step 1: Proving $w(q,p;0)\geq 0$ for $(q,p)\in[0,\pi]\times[p_{0},0]$.}

In this step, we choose a reflection parameter $\lambda\in(-\pi,0]$ and let the associated reflection function be
$$
w(q,p;\lambda)=h(q,p)-\tilde{h}(q,p)=h(q,p)-h(2\lambda-q,p),~~~for ~~~(q,p)\in[\lambda,2\lambda+\pi]\times[p_{0},0],
$$
it's obvious that the reflection function $w$ satisfies
\begin{eqnarray}\label{eq3.4}
\left\{\begin{array}{ll}
{w(\lambda,p;\lambda)=0,}& {\text { for } p\in[p_{0},0]}, \\
{w(q,p_{0};\lambda)=0,} & {\text { for } q\in[-\pi,\pi]}.\end{array}\right.
\end{eqnarray}
Since all streamlines are monotonic in a neighborhood of the trough line $q=-\pi$, there hold
$$
w(q,p;\lambda)\geq 0~~~for~~~|\lambda+\pi|<\varepsilon,
$$
where $\varepsilon$ is small enough, which means $\lambda$ is close enough to the trough line $q=-\pi$.
Thus, there exists an extremal position $\lambda_{0}$ for $\lambda$ by
$$
\lambda_{0}=\sup\{ \lambda\in(-\pi,0]: w(q,p;\lambda)\geq 0~~ for~~ all~~ (q,p)\in[\lambda,2\lambda+\pi]\times [p_{0},0] \}.
$$

\begin{figure}[ht]
\includegraphics[width=0.75\textwidth]{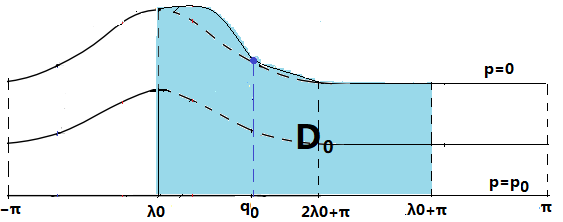}
\caption{The wave profile for $\lambda\in(-\pi,0]$}
\label{fig3.1}
\end{figure}

According to the assumption in the theorem, only one of the following cases can occur: (see Figure \ref{fig3.1} in the following)\\
(i) $\lambda_{0}=0$;\\
(ii) $\lambda_{0}\in(-\pi,0)$ and there exists $q_{0}\in [\lambda_{0},2\lambda_{0}+\pi]$ for which $w(q_{0},0;\lambda_{0})=0$, $w_{q}(q_{0},0;\lambda_{0})=0$ and $w_{qq}(q_{0},0;\lambda_{0})\geq0$.

Now let's preclude the second case (ii). For convenience, let's enlarge the domain of definition of $w$ by setting
$$
w(q,p;\lambda_{0})=h(q,p)-h(2\lambda_{0}+2\pi-q,p),~~for~~(q,p)\in[2\lambda_{0}+\pi,\lambda_{0}+\pi]\times [p_{0},0]
$$
and redefine
$$
D_{0}=(\lambda_{0},\lambda_{0}+\pi)\times(p_{0},0).
$$

At the same time, the periodicity of $h$ with respect to $q$ could guarantee that $w\in C^{2,\alpha}(\overline{D_{0}})$. In addition, if $2\lambda_{0}+\pi$ lies to
the left of the wave crest, the assumption of monotonicity for streamlines would prevent the second case (ii) occurring, which implies $\lambda_{0}=0$, that is to say $2\lambda_{0}+\pi=\pi$. This is contradicted with $2\lambda_{0}+\pi$ on the left of the wave crest. Thus, $2\lambda_{0}+\pi$ lies to the right of (or at least in line with) the wave crest. This observation provides that $h(q,p)$ is non-increasing for $q\in [2\lambda_{0}+\pi, \pi]$. Combining the fact and the definition of $\lambda_{0}$, we can obtain
\begin{equation}\label{eq3.5}
w(q,p;\lambda_{0})\geq 0 ~~~for~~~(q,p)\in D_{0}.
\end{equation}
In addition, we can also get the boundary condition
\begin{eqnarray}\label{eq3.6}
\left\{\begin{array}{ll}
{w(\lambda_{0},p;\lambda_{0})=w(\lambda_{0}+\pi,p;\lambda_{0})=0,} ~~~for~~~{ p\in[p_{0},0]}, \\
{w(q,p_{0};\lambda_{0})=0} ~~~and~~~ { w(q,0;\lambda_{0})\geq0, ~~~for~~~q\in[\lambda_{0},\lambda_{0}+\pi]}.\end{array}\right.
\end{eqnarray}
Based on (\ref{eq3.5}) and (\ref{eq3.6}), we can use the strong maximum principle in Lemma \ref{lem3.1} to show that the minimum of $w(q,p;\lambda_{0})$ must be attained on $\partial\Omega$ unless $w(q,p;\lambda_{0})$ is a constant, that is to say
\begin{equation}\label{eq3.7}
w(q,p;\lambda_{0})>0~~in~~D_{0}~~~or ~~~ w(q,p;\lambda_{0})\equiv 0~~in ~~\overline{D_{0}}.
\end{equation}

In the following, we will show that neither of (\ref{eq3.7}) can hold. Let's first assume by contradiction that $w(q,p;\lambda_{0})>0$ in $D_{0}$. Now we claim that
\begin{equation}\label{eq3.8}
w_{p}(q_{0},0;\lambda_{0})\geq 0.
\end{equation}
Indeed, at the point $(q_{0},0)$, we have $w(q_{0},0;\lambda_{0})=0$ and $w_{q}(q_{0},0;\lambda_{0})=0$, which implies that
\begin{equation}\label{eq3.9}
h(q_{0},0)=h(2\lambda_{0}-q_{0},0)=\tilde{h}(q_{0},0),~~~h_{q}(q_{0},0)=-h_{q}(2\lambda_{0}-q_{0},0)=-\tilde{h}_{q}(q_{0},0).
\end{equation}
Considering the second formula of (\ref{eq3.3}) and (\ref{eq3.9}), at $(q_{0},0)$, we have
\begin{equation}\label{eq3.10}
(2g\rho h-Q-2\sigma\frac{h_{qq}}{(1+h_{q}^{2})^{\frac{3}{2}}})h_{p}^{2}-(2g\rho h-Q-2\sigma\frac{\tilde{h}_{qq}}{(1+h_{q}^{2})^{\frac{3}{2}}})\tilde{h}_{p}^{2}=0.
\end{equation}
It follows that
\begin{equation}\label{eq3.11}
(2g\rho h-Q-2\sigma\frac{h_{qq}}{(1+h_{q}^{2})^{\frac{3}{2}}})(h_{p}^{2}-\tilde{h}_{p}^{2})=2\sigma\frac{h_{qq}-\tilde{h}_{qq}}{(1+h_{q}^{2})^{\frac{3}{2}}}
=2\sigma\frac{w_{qq}(q_{0},0;\lambda_{0})}{(1+h_{q}^{2})^{\frac{3}{2}}}\leq0,
\end{equation}
where the last inequality holds from $\sigma\leq0$ and $w_{qq}(q_{0},0;\lambda_{0})\geq 0$. Besides, on $p=0$, we also have that
$$
1+h_{q}^{2}+(2g\rho h-Q-2\sigma\frac{h_{qq}}{(1+h_{q}^{2})^{\frac{3}{2}}})h_{p}^{2}=0,
$$
that is to say,
\begin{equation}\label{eq3.12}
 2g\rho h-Q-2\sigma\frac{h_{qq}}{(1+h_{q}^{2})^{\frac{3}{2}}}=-\frac{1+h_{q}^{2}}{h_{p}^{2}}<0
\end{equation}
Then (\ref{eq3.11}) and (\ref{eq3.12}) yield that
\begin{equation}\label{eq3.13}
h_{p}^{2}-\tilde{h}_{p}^{2}=(h_{p}+\tilde{h}_{p})(h_{p}-\tilde{h}_{p})\geq0.
\end{equation}
From the fact $h_{p}=\frac{1}{\sqrt{\rho}(c-u)}>0$ and $\tilde{h}_{p}=\frac{1}{\sqrt{\rho}(c-\tilde{u})}>0$, then (\ref{eq3.8}) follows from
$$
h_{p}(q_{0},0)-\tilde{h}_{p}(q_{0},0)\geq 0.
$$
On the other hand, at the point $(q_{0},0)$, we use the Hopf lemma to obtain
$$
\partial_{\vec{\nu}}w|_{(q_{0},0)}<0,
$$
which is contradicted with the fact $w_{q}(q_{0},0;\lambda_{0})=0$ and (\ref{eq3.8}).

Therefore, in (\ref{eq3.7}), there maybe hold
$$
w(q,p;\lambda_{0})\equiv 0~~~in ~~~\overline{D_{0}}.
$$
This would imply that $h(2\lambda_{0}+\pi,0)=h(\pi,0)=h(-\pi,0)$, which contradicts with the assumption that the streamlines are locally strictly monotonic near the troughs. To sum up, only the first case (i) (i.e. $\lambda_{0}=0$) would occur. Then the definition of $\lambda_{0}$ gives that
$$
w(q,p;0)\geq 0~~~for ~~~(q,p)\in[0,\pi]\times[p_{0},0].
$$
~\\
{\bf Step 2: Proving $w(q,p;0)\leq 0$ for $(q,p)\in[0,\pi]\times[p_{0},0]$.}

Here we use the similar argument as {\bf Step 1}, for a reflection parameter, let's introduce the reflection function
$$
w(q,p;\lambda)=h(q,p)-\tilde{h}(q,p)=h(q,p)-h(2\lambda-q,p),~~~for ~~~(q,p)\in[2\lambda-\pi,\lambda]\times[p_{0},0],
$$
it's easy to see that
\begin{eqnarray}\label{eq3.14}
\left\{\begin{array}{ll}
{w(\lambda,p;\lambda)=0,}& {\text { for } p\in[p_{0},0]}, \\
{w(q,p_{0};\lambda)=0,} & {\text { for } q\in[-\pi,\pi]}.\end{array}\right.
\end{eqnarray}
Since all streamlines are monotonic in a neighborhood of the trough line $q=\pi$, there holds
$$
w(q,p;\lambda)\geq 0~~~for~~~|\lambda-\pi|<\varepsilon,
$$
where $\varepsilon$ is small enough, which means $\lambda$ is close enough to the trough line $q=\pi$.
Thus, there exists
$$
\lambda_{0}=\inf\{ \lambda\in[0,\pi): w(q,p;\lambda)\geq 0~~ for~~ all~~ (q,p)\in[2\lambda-\pi,\lambda]\times [p_{0},0] \}.
$$

\begin{figure}[ht]
\includegraphics[width=0.75\textwidth]{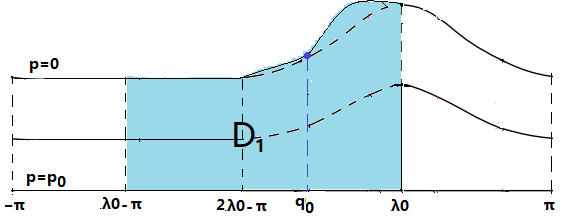}
\caption{The wave profile for $\lambda\in[0,\pi)$}
\label{fig3.2}
\end{figure}

One of the following two alternatives must occur: (see Figure \ref{fig3.2})\\

(i) $\lambda_{0}=0$;\\

(ii) $\lambda_{0}\in(0,\pi)$ and there exists $q_{0}\in [2\lambda_{0}-\pi,\lambda_{0}]$ for which $w(q_{0},0;\lambda_{0})=0$, $w_{q}(q_{0},0;\lambda_{0})=0$ and $w_{qq}(q_{0},0;\lambda_{0})\geq0$.

Now let's show the second case (ii) is impossible. As before, we first enlarge the domain of definition of $w$ by the following extension
$$
w(q,p;\lambda_{0})=h(q,p)-h(2\lambda_{0}-2\pi-q,p),~~~for~~~(q,p)\in [\lambda_{0}-\pi,2\lambda_{0}-\pi]\times [p_{0},0]
$$
and let
$$
D_{1}=(\lambda_{0}-\pi,\lambda_{0})\times (p_{0},0).
$$

Then it turns out that the periodicity of $h$ can ensure $w\in C^{2,\alpha}(\overline{D_{1}})$. If $2\lambda_{0}-\pi$ lies to the right of the wave crest, the assumption on the monotonicity of streamlines would imply that the first case (i) holds, which leads to the contradiction due to $2\lambda_{0}-\pi=-\pi$. Thus, $2\lambda_{0}-\pi$ lies to left of (or at least in line with) the wave crest. Consequently, we have that $h(q,p)$ is non-decreasing for $q\in[-\pi,2\lambda_{0}-\pi]$. Combining this fact and the definition of $\lambda_{0}$, we have
$$
w(q,p;\lambda_{0})\geq 0~~~for~~~(q,p)\in D_{1}.
$$
Similarly, we also have the boundary condition
\begin{eqnarray}\label{eq3.15}
\left\{\begin{array}{ll}
{w(\lambda_{0},p;\lambda_{0})=w(\lambda_{0}-\pi,p;\lambda_{0})=0,} ~~~for~~~{ p\in[p_{0},0]}, \\
{w(q,p_{0};\lambda_{0})=0} ~~~and~~~ { w(q,0;\lambda_{0})\geq0, ~~~for~~~q\in[\lambda_{0}-\pi,\lambda_{0}]}.\end{array}\right.
\end{eqnarray}
Then, using the maximum principle in in Lemma \ref{lem3.1}, we obtain
\begin{equation}\label{eq3.16}
w(q,p;\lambda_{0})>0~~in~~D_{1}~~~or ~~~ w(q,p;\lambda_{0})\equiv 0~~in ~~\overline{D_{1}}.
\end{equation}
If $w(q,p;\lambda_{0})>0$ in $D_{1}$, by using the Hopf lemma at point $(q_{0},0)$, we would get that
\begin{equation}\label{eq3.17}
\partial_{\vec{\nu}}w|_{(q_{0},0)}< 0
\end{equation}
On the other hand, we can also obtain that
\begin{equation}\label{eq3.18}
w_{p}(q_{0},0;\lambda_{0})\geq0
\end{equation}
by using the same argument in {\bf Step 1}. Then (\ref{eq3.17}) contradicts with the fact $w_{q}(q_{0},0;\lambda_{0})=0$ and (\ref{eq3.18}).

Thus, in (\ref{eq3.16}), there holds
$$
w(q,p;\lambda_{0})\equiv 0~~~in ~~~\overline{D_{1}},
$$
which gives that $h(-\pi,0)=h(2\lambda_{0}-\pi,0)=h(\pi,0)$, which is contradicted with our assumptions that streamlines are locally strictly monotonic near the troughs. Summing up, we obtain $\lambda_{0}=0$. Then the definition of $\lambda_{0}$ yields that
$$
w(q,p;0)\leq 0~~~for~~~(q,p)\in[0,\pi]\times[p_{0},0].
$$

\begin{remark}:
In particular, if choosing $\sigma=0$ in this case, the symmetry result is consistent with the recent work \cite{6}.
\end{remark}

\section{\bf The stratified water waves with stagnation points and $\sigma=0$ }

In this section, we don't consider the effect of surface tension and admit the internal stagnation points but choose $\rho$ to depend linearly on the streamline and the Benoulli function to be a constant. These choices are not merely a mathematical simplication. In fact, they are related to physical phenomena, which can be seen in \cite{14}. Then we will set
$$
\rho(-\psi)=A\psi+B
$$
and
$$
-\beta(\psi)=\gamma,
$$
where $A,\gamma\in R$ and $B\in R^{+}$ are constants. Indeed, the constant $B$ is also arbitrary but since the density function should be nonnegative all the times and $\psi=0$ at the surface, then we require that $B\in R^{+}$. Therefore, in this case, the governing equations (\ref{eq2.12}) would become
\begin{eqnarray}
\left\{\begin{array}{ll}{ \Delta\psi+Agy=\gamma} & {\text { in }0 < y<\eta(x)}, \\
{|\nabla\psi|^{2}+2gBy=Q} & { \text { on } y=\eta(x)}, \\
{\psi=0} & {\text { on } y=\eta(x)}, \\
{\psi=-p_{0}} & {\text { on } y=0},\end{array}\right. \label{eq4.1}
\end{eqnarray}
where $Q=2(E|_{\eta}-P_{atm})$ is a constant.

To extend the proof of symmetry when stagnation points occur, we first investigate the information on existence and location of stagnation points.

\subsection{The existence of stagnation points }
Consider the wave train with a flat surface and let $h$ be the height over the flat bed, $L$ be the wavelength, $k=\frac{2\pi}{L}$ be the wave number, $A,B,\gamma$ be the constants as in (\ref{eq4.1}) and $p_{0}$ is the pseudo mass flux defined in Section 2. It's known that the simplest solutions to the water wave problem (\ref{eq4.1}) are laminar flows as follows, which only depend on the vertical direction.
\begin{equation}\label{eq4.2}
\psi(x,y)=\frac{\gamma}{2}y^{2}-\frac{Ag}{6}y^{3}+\left(\frac{p_{0}}{h}-\frac{\gamma h}{2}+\frac{Agh^{2}}{6}\right)y-p_{0},~~for~~0\leq y\leq h,
\end{equation}
with the pseudo velocity field
\begin{equation}\label{eq4.3}
(\psi_{y},-\psi_{x})=\left(\gamma y-\frac{Ag}{2}y^{2}+\frac{p_{0}}{h}-\frac{\gamma h}{2}+\frac{Agh^{2}}{6},~0\right),~~for~~0\leq y\leq h.
\end{equation}
Furthermore, we can write (\ref{eq4.3}) as
\begin{equation}\label{eq4.4}
(\psi_{y},-\psi_{x})=\left((h-y)\left( \frac{Ag}{2}(y+h)-\gamma  \right)+\lambda,~0\right),~~for~~0\leq y\leq h,
\end{equation}
by letting $\lambda:=\frac{p_{0}}{h}+\frac{\gamma h}{2}-\frac{Agh^{2}}{3}$ (see \cite{14}).

Indeed, the laminar flows can give rise to genuine wave solutions to the water wave problem. It's known that determining the bifurcation values is significant by using the Crandall-Rabinowitz theorem to find nontrivial solutions. To attain this, we need to consider the solutions of the following equation
$$
k\coth(kh)\lambda^{2}-(\gamma-Agh)\lambda-gB=0,
$$
whose two solutions are
\begin{equation}\label{eq4.5}
\lambda_{\pm}=\frac{(\gamma-Agh)\tanh(kh)}{2k}\pm\sqrt{\frac{(\gamma-Agh)^{2}\tanh^{2}(kh)}{4k^{2}}+\frac{gB}{k}\tanh(kh)}.
\end{equation}
From (\ref{eq4.4}), we know that stagnation points exist in laminar flows if and only if
\begin{equation}\label{eq4.6}
(h-y)\left( \frac{Ag}{2}(y+h)-\gamma  \right)+\lambda=0.
\end{equation}
Therefore, from (\ref{eq4.5})-(\ref{eq4.6}), we obtain that
\begin{itemize}
  \item If $y=h$ or $\gamma=\frac{Ag}{2}(y+h)$, the existence of stagnation points is impossible due to $\lambda_{+}>0$ and $\lambda_{-}<0$. In fact, this implies that there are not stagnation points on the free surface of laminar flows.
  \item If $\gamma>\frac{Ag}{2}(y+h)$, it follows that stagnation points can only occur for $\lambda_{+}$ on the laminar flow inside the fluid.
  \item If $\gamma<\frac{Ag}{2}(y+h)$, it follows that stagnation points can only occur for $\lambda_{-}$ on the laminar flow inside the fluid.
\end{itemize}
By perturbation, it's obvious that if the wave amplitude is small there exist stagnation points on the streamline inside the fluid away from the free surface. Note that the existence and bifurcation structure of such waves can be ensured by bifurcation method in \cite{14}, which is based on the important work \cite{8,9,10}. From the physical view, it's reasonable to require that all the points with $u=c$ lie beneath the wave trough lines.

\subsection{ Symmetry of the stratified waves with stagnation points}
Without loss of generality, here we are seeking the one period $2\pi$ (i.e. $k=1$) of the wave-train for our governing system (\ref{eq4.1}). For simplicity, we choose the trough of the surface to be at $x=\pm\pi$ and we can obtain the following symmetry result.
\begin{theorem}\label{thm4.1}
The stream function $\psi\in C^{2}(\Omega)$ which is solution of the system (\ref{eq4.1}) is symmetric about $x=0$ if the wave profile is monotonic between the crest and the trough.
\end{theorem}

\begin{figure}[ht]
\includegraphics[width=0.75\textwidth]{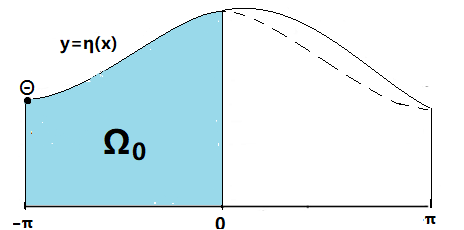}
\caption{The wave profile for $\lambda_{0}=0$}
\label{fig4.1}
\end{figure}

\begin{figure}[ht]
\includegraphics[width=0.75\textwidth]{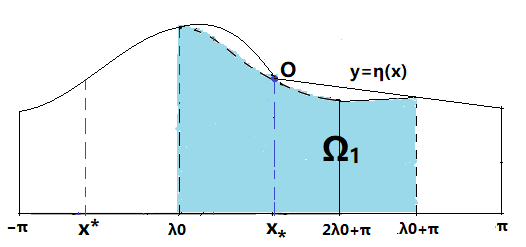}
\caption{The wave profile for $\lambda_{0}<0$}
\label{fig4.2}
\end{figure}

Before starting the proof, let's first introduce a suitable Serrin's Edge Point Lemma.
\begin{lemma} (cf. \cite{4})\label{lem4.1}
 Suppose $\Omega\subset R^{2}$ be the domain defined above and $L$ be an uniformly elliptic operator. Let $T$ be a line normal to the top boundary $\eta(x)$ at some point $\Theta$. Choose a part of $\Omega$ lying on a particular side of the line $T$ and denote it by $\Omega_{0}$. If there exist $w\in C^{2}(\overline{\Omega_{0}})$ with $Lw\leq 0$ and $w> 0$ in $\Omega_{0}$ and $w=0$ at $\Theta$, then
 $$
 either~~~~\partial_{\vec{\nu}}w(\Theta)<0~~~~or~~~~\partial_{\vec{\nu}}^{2}\omega(\Theta)<0,
 $$
 where $\vec{\nu}$ is a non-tangential outward vector at $\Theta$, unless $w\equiv0$ throughout $\overline{\Omega_{0}}$.
\end{lemma}
~\\
{\bf Proof of the Theorem \ref{thm4.1}:}
The similar ideas in \cite{15,16} would be adopted to finish our proof. Now let's first fix $\lambda\in(-\pi,0]$ and define
$$
\Omega_{\lambda}:=\{ (x,y)\in R^{2}|~-\pi<x<\lambda,~~0<y<\eta(x)  \}\subset\Omega
$$
and its reflected region about the line $x=\lambda$ by
$$
\Omega_{\lambda}^{r}:=\{ (2\lambda-x,y)\in R^{2}|~-\pi<x<\lambda,~~0<y<\eta(x)  \}.
$$
Since $x\pm \pi$ is the wave trough, the monotonicity assumption would ensure that the surface $\eta(x)$ is non-decreasing on the interval $(-\pi,-\pi+2\varepsilon)$ such that $\Omega_{\lambda}^{r}$ is a subset of fluid domain $\Omega$ for $\lambda\in (-\pi,-\pi+\varepsilon)$.
Now we move the line $x=\lambda$, there will be a critical value where the reflected domain is still a subset of $\Omega$. Define the critical value by
$$
\lambda_{0}=\max\{ \lambda|~ \Omega_{\lambda}^{r}\subset\Omega\}\leq 0
$$
and define
$$
\Omega_{0}=\{ (x,y)\in R^{2}|~-\pi<x<\lambda_{0},~~0<y<\eta(x)\}\subset\Omega.
$$
Similarly, there are two possible cases for $\lambda_{0}$:\\
(i) $\lambda_{0}=0$; (see Figure \ref{fig4.1})\\
(ii) $\lambda_{0}\in(-\pi,0)$ and there exist a $x_{*}\in[\lambda_{0},2\lambda_{0}+\pi]$ such that $\partial \Omega_{\lambda}^{r}$ is tangent to $\partial \Omega$ at point $(x_{*}, \eta(x_{*}))$. (see Figure \ref{fig4.2})

In the following, we first show that the wave crest is at $(0,\eta(0))$ and $\psi(x,y)$ is symmetric about $x=0$ if the case (i) occurs. Then we preclude the possibility of case (ii) by contradiction analysis.\\
~\\
{\bf Case (i):}

In this case, we have $\lambda_{0}=0$. Based on the analysis in subsection 4.1, it's known that $\psi_{y}|_{\eta(x)}\neq 0$. Without loss of generality, we assume $\psi_{y}|_{\eta(x)}> 0$ and introduce the function $m(x,y)$ by
$$
m(x,y)=\psi(x,y)-\psi(-x,y)
$$
on the domain
$$
\Omega_{0}:=\{ (x,y)\in R^{2}|~-\pi<x<0,~ 0<y<\eta(x) \}.
$$

It's easy to see that $m\in C^{2}(\overline{\Omega_{0}})$ and there hold
\begin{equation}\label{eq4.7}
\Delta m(x,y)=\Delta\psi(x,y)-\Delta\psi(-x,y)=0, ~~in~~\Omega_{0}
\end{equation}
due to (\ref{eq4.1}). Now let's pay attention to the situation on the boundary. It's obvious that
\begin{eqnarray}
\left\{\begin{array}{ll}{ m(x,0)=-p_{0}-(-p_{0})=0}, \\
{m(-\pi,y)=\psi(-\pi,y)-\psi(\pi,y)=0}, \\
{m(0,y)=\psi(0,y)-\psi(0,y)=0}.\end{array}\right. \label{eq4.8}
\end{eqnarray}
Then we look at the values of $m$ on the surface
\begin{equation}\label{eq4.9}
m(x,\eta(x))=\psi(x,\eta(x))-\psi(-x,\eta(x))=0-\psi(-x,\eta(x))\geq 0
\end{equation}
due to the definition of $\lambda_{0}$, $\psi_{y}|_{\eta(x)}>0$ and its continuity.
From (\ref{eq4.7})-(\ref{eq4.9}), we can apply the strong maximum principle to obtain that
\begin{equation}\label{eq4.10}
m>0~~in~~\Omega_{0}
\end{equation}
unless
\begin{equation}\label{eq4.11}
m\equiv 0~~on~~\overline{\Omega_{0}}.
\end{equation}

In the following, we will use Lemma \ref{lem4.1} to establish the latter case. Let's choose $\Theta:=(-\pi,\eta(-\pi))$ and $T:=\{ x=-\pi \}$. It's easy to see that the line $T$ is normal to the surface $\eta(x)$ at point $\Theta$. If (\ref{eq4.10}) holds, we will compute all partial derivatives of $m$ up to order two at $\Theta$ to deduce the contradiction.

By the definition of $m(x,y)$, it follows that
\begin{equation}\label{eq4.12}
m(\Theta)=0,~~m_{y}(\Theta)=0,~~m_{yy}(\Theta)=0.
\end{equation}
Differentiating $\psi(x,\eta(x))=0$ with respect to $x$, we obtain that
$$
\psi_{x}+\psi_{y}\eta_{x}(x)=0.
$$
Due to $\eta_{x}(x)=0$ at wave trough $\Theta$, then we have that
$$
\psi_{x}(\Theta)=0,
$$
which yields that
\begin{equation}\label{eq4.13}
m_{x}(\Theta)=2\psi_{x}(\Theta)=0.
\end{equation}
Due to the $2\pi$-periodicity of $\psi$ about $x$, we have that
\begin{equation}\label{eq4.14}
m_{xx}(\Theta)=[\psi_{xx}(x,y)-\psi_{xx}(-x,y)]|_{\Theta}=0.
\end{equation}
As last, let's differentiate the nonlinear boundary condition $|\nabla\psi|^{2}+2gBy=Q$ with respect to $x$ and evaluate the result at point $\Theta$, which gives that
$$
[2\psi_{x}(\psi_{xx}+\psi_{xy}\eta_{x})+2\psi_{y}(\psi_{xy}+\psi_{yy}\eta_{x})]|_{\Theta}=0,
$$
that is to say,
$$
2\psi_{y}(\Theta)\psi_{xy}(\Theta)=0.
$$
Since $\psi_{y}(\Theta)>0$, then
$$
\psi_{xy}(\Theta)=0,
$$
which yields
\begin{equation}\label{eq4.15}
m_{xy}(\Theta)=2\psi_{xy}(\Theta)=0.
\end{equation}
Combining (\ref{eq4.12})-(\ref{eq4.15}) with Serrin's Edge Point lemma, we deduce that only (\ref{eq4.11}) holds, hence symmetry is attained.
By the way, if $\psi_{y}|_{\eta(x)}<0$, the same process can be repeated by defining in $\Omega_{0}$
$$
m(x,y)=\psi(-x,y)-\psi(x,y).
$$

{\bf Case (ii):}

In this case, the reflected surface is tangent to $\partial\Omega$ at point $O:=(x_{*},\eta(x_{*}))$. Similarly, assume $\psi_{y}|_{\eta(x)}>0$, now we redefine $m(x,y)$ by
\begin{eqnarray}
m(x,y)=\left\{\begin{array}{ll}
{\psi(2\lambda_{0}-x,y)-\psi(x,y),}& {\text { for }\lambda_{0}\leq x\leq 2\lambda_{0}+\pi,~0\leq y\leq \eta(2\lambda_{0}-x)}, \nonumber \\
{\psi(2\lambda_{0}+2\pi-x,y)-\psi(x,y),}& {\text { for }2\lambda_{0}+\pi\leq x\leq \lambda_{0}+\pi,~0\leq y\leq \eta(2\lambda_{0}+2\pi-x)}.\end{array}\right.  \nonumber
\end{eqnarray}
and define the domain
$$
\Omega_{1}:=\{(x,y)\in(\lambda_{0},2\lambda_{0}+\pi]\times(0,\eta(2\lambda_{0}-x))\}\cup\{(x,y)\in(2\lambda_{0}+\pi,\lambda_{0}+\pi)\times(0,\eta(2\lambda_{0}+2\pi-x))\}
$$

In fact, here we have to deal with two reflections. As before, it's obvious that
\begin{eqnarray}\label{eq4.16}
\left\{\begin{array}{llll}{\Delta m=0} & {\text { in } \Omega_{1} }, \\
{m(x,0)=0} & {\text { for } \lambda_{0} \leq x \leq \lambda_{0}+\pi}, \\
{m(\lambda_{0},y)=0} & {\text { for }0 \leq y\leq \eta(\lambda_{0})}, \\
{m(\lambda_{0}+\pi,y)=0} & { \text { for }0 \leq y \leq \eta(\lambda_{0}+\pi)}, \\
{m(x,\eta(2\lambda_{0}-x))\geq 0} & {\text { for } \lambda_{0}\leq x\leq 2\lambda_{0}+\pi}, \\
{m(x,\eta(2\lambda_{0}+2\pi-x))\geq 0} & {\text { for } 2\lambda_{0}+\pi\leq x\leq \lambda_{0}+\pi}\end{array}\right.
\end{eqnarray}
due to the definition of $\lambda_{0}$, $\psi_{y}|_{\eta(x)}>0$ and its continuity. Based on (\ref{eq4.16}), we use the strong maximum principle again to get that
\begin{equation}\label{eq4.17}
m>0~~in~~\Omega_{1}
\end{equation}
unless
\begin{equation}\label{eq4.18}
m\equiv 0~~on~~\overline{\Omega_{1}}.
\end{equation}

In the following, we will establish that, if (\ref{eq4.17}) occurs, there will be a contradiction at point $O$ by using Hopf lemma. If (\ref{eq4.18}) holds, which is contradicted with monotonicity assumption. Thus, we preclude the possibility of {\bf Case (ii)}.

Now let's show that how to use Hopf lemma at point $O$. It's known that $O=(x_{*},\eta(x_{*}))$ is the tangency point for $x_{*}\in [\lambda_{0},2\lambda_{0}+\pi]$. It's easy to see that
\begin{equation}\label{eq4.19}
m(O)=m(x_{*},\eta(x_{*}))=\psi(2\lambda_{0}-x_{*},\eta(x_{*}))-\psi(x_{*},\eta(x_{*})=0
\end{equation}
because both points are on the surface. From (\ref{eq4.17}) and (\ref{eq4.19}), we use the Hopf lemma to obtain that
\begin{equation}\label{eq4.20}
\partial m_{\vec{\nu}}(O)<0.
\end{equation}

On the other hand, if let $x^{*}=2\lambda_{0}-x_{*}$, we have that
\begin{equation}\label{eq4.21}
\eta(x_{*})=\eta(x^{*})~~~and~~~\eta_{x}(x_{*})=-\eta_{x}(x^{*})
\end{equation}
due to the particularity of $O$ (see Figure \ref{fig4.2}). Considering the nonlinear boundary condition at $x=x_{*}$ and $x=x^{*}$, we can obtain that
\begin{equation}\label{eq4.22}
\psi_{x}^{2}(x_{*},\eta(x_{*}))+\psi_{y}^{2}(x_{*},\eta(x_{*}))=\psi_{x}^{2}(x^{*},\eta(x^{*}))+\psi_{y}^{2}(x^{*},\eta(x^{*}))
\end{equation}
due to (\ref{eq4.21}). In addition, if differentiating $\psi(x,\eta(x))=0$ with respect to $x$ and evaluating at $x=x_{*}$ and $x=x^{*}$, we have that
\begin{equation}\label{eq4.23}
\frac{\psi_{x}(x_{*},\eta(x_{*}))}{\psi_{y}(x_{*},\eta(x_{*}))}=-\eta_{x}(x_{*})=\eta_{x}(x^{*})=-\frac{\psi_{x}(x^{*},\eta(x^{*}))}{\psi_{y}(x^{*},\eta(x^{*}))},
\end{equation}
where (\ref{eq4.21}) is used again. Combining (\ref{eq4.22}) and (\ref{eq4.23}) with the fact $\psi_{y}|_{\eta(x)}>0$, we obtain that
\begin{equation}\label{eq4.24}
\psi_{x}(x_{*},\eta(x_{*}))=-\psi_{x}(x^{*},\eta(x^{*})),~~ ~~\psi_{y}(x_{*},\eta(x_{*}))=\psi_{y}(x^{*},\eta(x^{*})).
\end{equation}
It's obvious that (\ref{eq4.21}) and (\ref{eq4.24}) imply
$$
m_{x}(O)=-\psi_{x}(x^{*},\eta(x_{*}))-\psi_{x}(x_{*},\eta(x_{*}))=-\psi_{x}(x^{*},\eta(x^{*}))-\psi_{x}(x_{*},\eta(x_{*}))=0
$$
and
$$
m_{y}(O)=\psi_{y}(x^{*},\eta(x_{*}))-\psi_{y}(x_{*},\eta(x_{*}))=\psi_{y}(x^{*},\eta(x^{*}))-\psi_{y}(x_{*},\eta(x_{*}))=0,
$$
which is contradicted with (\ref{eq4.20}). By the way, if $\psi_{y}|_{\eta(x)}<0$, the similar process can be carried out by letting
\begin{eqnarray}
m(x,y)=\left\{\begin{array}{ll}
{\psi(x,y)-\psi(2\lambda_{0}-x,y),}& {\text { for }\lambda_{0}\leq x\leq 2\lambda_{0}+\pi,~0\leq y\leq \eta(2\lambda_{0}-x)}, \nonumber \\
{\psi(x,y)-\psi(2\lambda_{0}+2\pi-x,y),}& {\text { for }2\lambda_{0}+\pi\leq x\leq \lambda_{0}+\pi,~0\leq y\leq \eta(2\lambda_{0}+2\pi-x)}.\end{array}\right.  \nonumber
\end{eqnarray}
~\\

\begin{remark}:
The Theorem \ref{thm4.1} and the fact stagnation points away from surface imply that the wave profile described by the function $\eta(x)$ is symmetric between the crest and the trough. Indeed, from (\ref{eq2.4}) and (\ref{eq2.5}), we have that
$$
v=(u-c)\eta_{x},~~\psi_{x}=-\sqrt{\rho}v,~~~\psi_{y}=\sqrt{\rho}(u-c),~~~on~~y=\eta(x).
$$
 Thus we obtain that
$$
\eta_{x}=\frac{v}{u-c}=\frac{-\psi_{x}}{\psi_{y}},~~~on~~y=\eta(x)
$$
is odd about $x$ because $\psi(x,y)$ is even about $x$, which yields that $\eta(x)$ is even about $x$.
In particular, if assuming $A=0$ in this section, the symmetry result is consistent with the case \cite{16}.
\end{remark}

\section*{Acknowledgments}
This work was partially supported by the National Natural Science Foundation of China (No.11571057).

\section*{Compliance with ethical standards}
{\bf Conflict of interest} The authors declare that they have no conflict of interest.

\end{document}